\documentclass[12pt]{amsart}%
\usepackage{amsmath}
\usepackage{graphicx}
\usepackage{amsfonts}
\usepackage{amssymb}%
\newtheorem{theorem}{Theorem}
\theoremstyle{plain}

\newtheorem{lemma}{Lemma}

\numberwithin{equation}{section}
\begin{document}
\title[Keller-Segel instability]{Pattern formation (I): The Keller-Segel Model}
\author{Yan Guo}
\address{Division of Applied Mathematics, Brown University, Providence, RI 02912, USA}
\email{guoy@dam.brown.edu}
\author{Hyung Ju Hwang}
\address{School of Mathematics, Trinity College, Dublin 2, Ireland}
\email{hjhwang@maths.tcd.ie}
\date{}

\begin{abstract}
We investigate nonlinear dynamics near an unstable constant equilibrium in the
classical Keller-Segel model. Given any general perturbation of magnitude
$\delta$, we prove that its nonlinear evolution is dominated by the
corresponding linear dynamics along a fixed finite number of fastest growing
modes, over a time period of $\ln\frac{1}{\delta}$. Our result can be
interpreted as a rigourous mathematical characterization for early pattern
formation in the Keller-Segel model.

\end{abstract}
\maketitle

\section{ Growing Modes in the Keller-Segel Model}

The goal of this section is to review the well-known instability criterion for
the classical Keller-Segel model, which describes directed movement of
microorganisms and cells stimulated by the chemical which they produce
themselves. The Keller-Segel system takes the form%
\begin{align}
U_{t}  &  =-\nabla\left(  -\mu\nabla U+\chi U\nabla V\right)  ,\label{KS}\\
V_{t}  &  =\nabla\left(  D\nabla V\right)  +fU-kV,\nonumber
\end{align}
where $U(x,t)$ is the cell density, $V(x,t)$ the chemo-attractant, $\mu>0$ the
amoeboid motility, $\chi>0$ the chemotactic sensitivity, $D>0$ the diffusion
rate of cAMP, $f>0$ the rate of cAMP secretion per unit density of amoebae,
$k>0$ the rate of degradation of cAMP in environment.

We assume Neumann boundary conditions for $U(x,t)$ and $V(x,t)$, in a
$d$-dimensional box $x\in%
{\mathbb T}^d%
=(0,\pi)^{d},$ $d=1,2,3,$ i.e.,%
\begin{equation}
\frac{\partial U}{\partial x_{i}}=\frac{\partial V}{\partial x_{i}}=0,\text{
\ at }x_{i}=0,\pi\text{, for }1\leq i\leq d. \label{Neumann}%
\end{equation}
A uniform constant solution
\[
U(x,t)\equiv\bar{U},\;V(x,t)\equiv\bar{V}%
\]
forms a homogeneous steady state provided
\begin{equation}
f\bar{U}=k\bar{V}. \label{steady-state}%
\end{equation}

In this article, we study the nonlinear evolution of a perturbation
\[
u(x,t)=U(x,t)-\bar{U},\;v(x,t)=V(x,t)-\bar{V}%
\]
around $[\bar{U},\bar{V}]$, which satisfies the equivalent Keller-Segel
system:%
\begin{align}
u_{t} &  =\mu\nabla^{2}u-\chi\bar{U}\nabla^{2}v-\chi\nabla(u\nabla
v),\label{nonlinear1}\\
v_{t} &  =D\nabla^{2}v+fu-kv.\label{nonlinear2}%
\end{align}
The corresponding linearized Keller-Segel system then takes the form%
\begin{align}
u_{t} &  =\mu\nabla^{2}u-\chi\bar{U}\nabla^{2}v,\label{linear1}\\
v_{t} &  =D\nabla^{2}v+fu-kv.\label{linear2}%
\end{align}
We use $[\cdot,\cdot]$ to denote a column vector, and let
\[
\mathbf{w}(x,t)\equiv\lbrack u(x,t),v(x,t)].
\]
Let $\mathbf{q=}\left(  q_{1},..,q_{d}\right)  \in\Omega=\left(
{\mathbb N}%
\cup\left\{  0\right\}  \right)  ^{d}$ and let%
\[
e_{\mathbf{q}}(x)\equiv\prod_{i=1}^{d}\cos\left(  q_{i}x_{i}\right)  ,
\]
Then $\left\{  e_{\mathbf{q}}(x)\right\}  _{\mathbf{q\in}\Omega}$ forms a
basis of the space of functions in $%
{\mathbb T}^d%
$ that satisfy Neumann boundary conditions (\ref{Neumann}). We look for a
normal mode to the linear Keller-Segel system (\ref{linear1}) and
(\ref{linear2}) of the following form:%
\begin{equation}
\mathbf{w}\left(  x,t\right)  =\mathbf{r}_{\mathbf{q}}\exp\left(
\lambda_{\mathbf{q}}t\right)  e_{\mathbf{q}}(x),\label{normal-mode}%
\end{equation}
where $\mathbf{r}_{\mathbf{q}}$ is a vector depending on $\mathbf{q.}$
Plugging (\ref{normal-mode}) into (\ref{linear1})-(\ref{linear2}) yields%
\[
\lambda_{\mathbf{q}}\mathbf{r}_{\mathbf{q}}=\left(
\begin{array}
[c]{cc}%
-\mu q^{2} & \chi\bar{U}q^{2}\\
f & -Dq^{2}-k
\end{array}
\right)  \mathbf{r}_{\mathbf{q}},
\]
where $q^{2}=\sum_{i=1}^{d}q_{i}^{2}$. A nontrivial normal mode can be
obtained by setting%
\[
\det\left(
\begin{array}
[c]{cc}%
\lambda_{\mathbf{q}}+\mu q^{2} & -\chi\bar{U}q^{2}\\
-f & \lambda_{\mathbf{q}}+Dq^{2}+k
\end{array}
\right)  =0.
\]
This leads to the following dispersion formula for $\lambda_{\mathbf{q}}$:%
\begin{equation}
\lambda_{\mathbf{q}}^{2}+\{q^{2}\left(  \mu+D\right)  +k\}\lambda_{\mathbf{q}%
}+q^{2}\{\mu\left(  Dq^{2}+k\right)  -\chi\bar{U}f\}=0.\label{dispersion}%
\end{equation}
Thus we deduce the following well-known aggregation (i.e., linear instability)
criterion by requiring there exists a $q$ such that%
\begin{equation}
\mu\left(  Dq^{2}+k\right)  -\chi\bar{U}f<0,\label{instability-criterion}%
\end{equation}
to ensure that (\ref{dispersion}) has at least one positive root
$\lambda_{\mathbf{q}}.$ This clearly implies that $\mu k-\chi\bar{U}f<0,$ and
an elementary computation of the discriminant yields:%
\begin{align*}
&  \{q^{2}\left(  \mu+D\right)  +k\}^{2}-4q^{2}\{\mu\left(  Dq^{2}+k\right)
-\chi\bar{U}f\}\\
&  =q^{4}\left(  \mu-D\right)  ^{2}+k^{2}+2q^{2}\left(  \mu+D\right)
k+4q^{2}\{-\mu k+\chi\bar{U}f\}\\
&  >0
\end{align*}
for $q.$ Therefore, there exist two distinct real roots for all $\mathbf{q}$
to the quadratic equation (\ref{dispersion}), which we denote
\[
\lambda_{-}(\mathbf{q})<\lambda_{+}(\mathbf{q}).
\]
We denote the corresponding (linearly independent) eigenvectors by
$\mathbf{r}_{-}(\mathbf{q})$ and $\mathbf{r}_{+}(\mathbf{q}),$ such that
\begin{equation}
\mathbf{r}_{\pm}(\mathbf{q})=\left[  \frac{\lambda_{\pm}(\mathbf{q}%
)+D\mathbf{q}^{2}+k}{f},1\right]  .\label{r}%
\end{equation}
Clearly, for $q$ large,
\[
\mu\left(  Dq^{2}+k\right)  -\chi\bar{U}f>0.
\]
Hence there are only finitely many $\mathbf{q}$ such that $\lambda
_{+}(\mathbf{q})>0.$ We therefore denote the largest eigenvalue by
$\lambda_{\text{max}}>0$ and define
\[
\Omega_{\text{max}}\equiv\{\mathbf{q}\in\Omega\text{ such that }\lambda
_{+}(\mathbf{q})=\lambda_{\text{max }}\}.
\]
It is easy to see that there is one $q^{2}$ (possibly two) having
$\lambda_{\mathbf{q}}^{+}\left(  q^{2}\right)  =\lambda_{\text{max}}$ when we
regard $\lambda_{\mathbf{q}}^{+}$ as a function of $q^{2}$. We also denote
$\nu>0$ to be the gap between the $\lambda_{\text{max}}$ and the rest.

Given any initial perturbation $\mathbf{w}\left(  \mathbf{x},0\right)  $, we
can expand it as%
\[
\mathbf{w}\left(  \mathbf{x},0\right)  =\sum_{_{\mathbf{q}}\in\Omega
}\mathbf{w}_{\mathbf{q}}e_{\mathbf{q}}(x)=\sum_{_{\mathbf{q}}\in\Omega
}\{w_{\mathbf{q}}^{-}\mathbf{r}_{-}(\mathbf{q})+w_{\mathbf{q}}^{+}%
\mathbf{r}_{+}(\mathbf{q})\}e_{\mathbf{q}}(x),
\]
so that
\begin{equation}
\mathbf{w}_{\mathbf{q}}=w_{\mathbf{q}}^{-}\mathbf{r}_{-}(\mathbf{q}%
)+w_{\mathbf{q}}^{+}\mathbf{r}_{+}(\mathbf{q}).\label{inde}%
\end{equation}
The unique solution $\mathbf{w}\left(  x,t\right)  =[u\left(  x,t\right)
,v\left(  x,t\right)  ]$ to (\ref{linear1})-(\ref{linear2}) is given by
\begin{align}
\mathbf{w}\left(  x,t\right)   &  =\sum_{_{\mathbf{q}}\in\Omega}%
\{w_{\mathbf{q}}^{-}\mathbf{r}_{-}(\mathbf{q})\exp\left(  \lambda_{\mathbf{q}%
}^{-}t\right)  +w_{\mathbf{q}}^{+}\mathbf{r}_{+}(\mathbf{q})\exp\left(
\lambda_{\mathbf{q}}^{+}t\right)  \}e_{\mathbf{q}}(x)\label{l}\\
&  \equiv e^{\mathcal{L}t}\mathbf{w}\left(  x,0\right)  .\nonumber
\end{align}
For any $\mathbf{u}\left(  \cdot\mathbf{,}t\right)  \in\left[  L^{2}\left(
{\mathbb T}^d%
\right)  \right]  ^{2}$, we denote $\left\Vert \mathbf{u}\left(
\cdot\mathbf{,}t\right)  \right\Vert \equiv\left\Vert \mathbf{u}\left(
\cdot\mathbf{,}t\right)  \right\Vert _{L^{2}}$. Our main result of this
section is

\begin{lemma}
\label{lineargrowth}Assume the instability criterion
(\ref{instability-criterion}) is valid. Suppose
\[
\mathbf{w}\left(  x,t\right)  =[u\left(  x\mathbf{,}t\right)  ,v\left(
x\mathbf{,}t\right)  ]\equiv e^{\mathcal{L}t}\mathbf{w}\left(  x,0\right)
\]
as in (\ref{l}) is a solution to the linearized KS system (\ref{linear1}%
)-(\ref{linear2}) with initial condition $\mathbf{w}\left(  \mathbf{x}%
,0\right)  $. Then there exists a constant $C_{1}\geq1$ depending on
$k,\bar{U},D,\mu,f,\chi,$ such that%
\[
\left\Vert \mathbf{w}\left(  \mathbf{\cdot},t\right)  \right\Vert \leq
C_{1}\exp\left(  \lambda_{\text{max}}t\right)  \left\Vert \mathbf{w}\left(
\mathbf{\cdot},0\right)  \right\Vert ,
\]
for all $t\geq0$.
\end{lemma}

\begin{proof}
We first consider the case for $t\geq1.$ By analyzing (\ref{dispersion}), for
$q$ large, we have
\[
\lim_{q\rightarrow\infty}\frac{\lambda_{\mathbf{q}}^{\pm}}{q^{2}}=-\mu,-D
\]
respectively. Notice that from the quadratic formula for (\ref{dispersion}),
\[
\frac{\lambda_{\mathbf{q}}^{+}-\lambda_{\mathbf{q}}^{-}}{q^{2}}\geq
\frac{2\sqrt{-\mu k+\chi\bar{U}f}}{q}.
\]
From solving (\ref{inde})
\begin{align*}
|w_{\mathbf{q}}^{\pm}| &  \leq\frac{1}{\det[\mathbf{r}_{-}(\mathbf{q}%
),\mathbf{r}_{+}(\mathbf{q})]}|\mathbf{r}_{\pm}(\mathbf{q})|\times
|\mathbf{w}_{\mathbf{q}}|\\
&  \leq Cq|\mathbf{w}_{\mathbf{q}}|,
\end{align*}
we deduce that for $t\geq1$ and $q$ large$,$%
\[
|w_{\mathbf{q}}^{\pm}\mathbf{r}_{\pm}(\mathbf{q})\exp\left(  \lambda
_{\mathbf{q}}^{\pm}t\right)  |\leq Cq|\mathbf{w}_{\mathbf{q}}|\exp(-\min
\{\mu,D\}q^{2}t)\leq C|\mathbf{w}_{\mathbf{q}}|.
\]
Thus we deduce the Lemma on the linear growth rate for $t\geq1$ by the formula
(\ref{l}).

On the other hand, for finite\ time $t\leq1,$ it suffices to derive the
standard energy estimate in $L^{2}$. From the Neumann boundary conditions, we
can take $u\times$ (\ref{linear1}) and add $Av\times$ of (\ref{linear2}) to
get%
\begin{align*}
&  \frac{1}{2}\frac{d}{dt}\int_{%
{\mathbb T}^d%
}\left\{  |u|^{2}+A|v|^{2}\right\} \\
&  +\int_{%
{\mathbb T}^d%
}\{\mu\left\vert \nabla u\right\vert ^{2}+AD\left\vert \nabla v\right\vert
^{2}\mathbf{-}\chi\bar{U}\nabla v\nabla u\}+Ak\int_{%
{\mathbb T}^d%
}|v|^{2}\\
&  =\int_{%
{\mathbb T}^d%
}Afuv.
\end{align*}
The integrand of the second integral can be chosen non-negative%
\begin{equation}
\mu\left\vert \nabla u\right\vert ^{2}+AD\left\vert \nabla v\right\vert
^{2}\mathbf{-}\chi\bar{U}\nabla v\nabla u\geq\frac{\mu}{2}\left\vert \nabla
u\right\vert ^{2}+\frac{\left(  \bar{U}\chi\right)  ^{2}\left\vert \nabla
v\right\vert ^{2}}{2\mu}\geq0, \label{lowerbound}%
\end{equation}
if the constant $A$ is
\begin{equation}
A=\frac{\left(  \bar{U}\chi\right)  ^{2}}{D\mu}. \label{a}%
\end{equation}
It thus follows that%
\[
\frac{1}{2}\frac{d}{dt}\int_{%
{\mathbb T}^d%
}\left\{  |u|^{2}+A|v|^{2}\right\}  \leq\frac{Af}{2}\int_{%
{\mathbb T}^d%
}\left\{  |u|^{2}+|v|^{2}\right\}  ,
\]
and the Gronwall inequality implies
\[
\left\Vert \mathbf{w}\left(  \mathbf{\cdot},t\right)  \right\Vert \leq
C\exp\left(  Ct\right)  ||\mathbf{w}(\mathbf{\cdot},0)||,
\]
for some $C>0.$ This immediately implies our lemma when $t\leq1.$
\end{proof}

\section{\bigskip Main Result}

Let $\theta$ be a small fixed constant, and $\lambda_{\text{max}}$ be the
dominant eigenvalue which is the maximal growth rate.\ We also denote the gap
between the largest growth rate $\lambda_{\text{max}}$ and the rest by
$\nu>0.$ Then for $\delta>0$ arbitrary small, we define the escape time
$T^{\delta}$ by%
\begin{equation}
\theta=\delta\exp\left(  \lambda_{\text{max}}T^{\delta}\right)  ,
\label{theta}%
\end{equation}
or equivalently%
\[
T^{\delta}=\frac{1}{\lambda_{\text{max}}}\ln\frac{\theta}{\delta}.
\]

Our main theorem is

\begin{theorem}
Assume that the set of $q^{2}=\sum_{i=1}^{d}q_{i}^{2}$ satisfying instability
criterion (\ref{instability-criterion}) is not empty for given parameters
$\mu,D,k,\chi,f$ and $\bar{U}$. Let%
\[
\mathbf{w}_{0}(x)=\sum_{_{\mathbf{q}}\in\Omega}\{w_{\mathbf{q}}^{-}%
\mathbf{r}_{-}(\mathbf{q})+w_{\mathbf{q}}^{+}\mathbf{r}_{+}(\mathbf{q}%
)\}e_{\mathbf{q}}(x).
\]
$\in H^{2}$ such that $||\mathbf{w}_{0}||=1.$ Then there exist constants
$\delta_{0}>0,$ $C>0,$ and $\theta>0,$ depending on $k,\bar{U},D,\mu,f,\chi,$
such that for all $0<\delta\leq\delta_{0}$, if the initial perturbation of the
steady state $[\bar{U},\bar{V}]$ in (\ref{steady-state}) is
\[
\mathbf{w}^{\delta}\left(  \mathbf{x},0\right)  =\delta\mathbf{w}_{0},
\]
then its nonlinear evolution $\mathbf{w}^{\delta}(t,x)$ satisfies
\begin{align*}
&  ||\mathbf{w}^{\delta}(t,x)-\delta e^{\lambda_{\max}t}\sum_{_{\mathbf{q}}%
\in\Omega_{\text{max}}}w_{\mathbf{q}}^{+}\mathbf{r}_{+}(\mathbf{q}%
)e_{\mathbf{q}}(x)||\\
&  \leq C\{e^{-\nu t}+\delta||w_{0}||_{H^{2}}^{2}+\delta e^{\lambda_{\max}%
t}\}\delta e^{\lambda_{\max}t}%
\end{align*}
for $0\leq t\leq T^{\delta},$ and $\nu>0$ is the gap between $\lambda_{\max}$
and the rest of $\lambda_{\mathbf{q}}$ in (\ref{dispersion}).
\end{theorem}

We notice that for $0\leq t\leq T^{\delta},$ $\delta e^{\lambda_{\max}t}%
\leq\theta,\,$\ is$\ $sufficiently small. As long as $w_{\mathbf{q}_{0}}%
^{+}\neq0$ for at least one $\mathbf{q}_{0}\mathbf{\in}\Omega_{\text{max}},$
which is generic for perturbations, the corresponding fastest growing modes
\[
||\delta e^{\lambda_{\max}t}\sum_{_{\mathbf{q}}\in\Omega_{\text{max}}%
}w_{\mathbf{q}}^{+}\mathbf{r}_{+}(\mathbf{q})e_{\mathbf{q}}||\geq\delta
e^{\lambda_{\max}t}|w_{\mathbf{q}_{0}}^{+}||\mathbf{r}_{+}(\mathbf{q}_{0})|,
\]
have the dominant leading order of $\delta e^{\lambda_{\max}t}.$ Our theorem
implies that the dynamics of a general perturbation is characterized by such
linear dynamics over a long time period of $\varepsilon T^{\delta}\leq t\leq
T^{\delta},$ for any $\varepsilon>0$. In particular, choose a fixed
$\mathbf{q}_{0}\in\Omega_{\text{max}}$ and let
\[
w_{0}(x)=\frac{\mathbf{r}_{+}(\mathbf{q}_{0})}{|\mathbf{r}_{+}(\mathbf{q}%
_{0})|}e_{\mathbf{q}_{0}}(x)
\]
then if $t=T^{\delta},$
\[
\left\Vert \mathbf{w}^{\delta}(t,\cdot)-\delta e^{\lambda_{\max}T^{\delta}%
}\frac{\mathbf{r}_{+}(\mathbf{q}_{0})}{|\mathbf{r}_{+}(\mathbf{q}_{0}%
)|}e_{\mathbf{q}_{0}}(\cdot)\right\Vert \leq C\{\delta^{\nu/\lambda_{\text{max
}}}+\theta^{2}\},
\]
hence
\[
\left\Vert \mathbf{w}^{\delta}(t,\cdot)\right\Vert \geq\theta-C\{\delta
^{\nu/\lambda_{\text{max }}}+\theta^{2}\}\geq\theta/2>0,
\]
which implies nonlinear instability as $\delta\rightarrow0$. The instability
occurs before the possible blow-up time.

In the early work of Keller and Segel \cite{KS} in 1970, they formulated the
advection-diffusion system (\ref{KS}) which consists of two parabolic
equations and viewed the initiation of Slime mold aggregation as instability.
Linearized system was used to analyze early stage of pattern formation and its
instability around homogeneous steady states. This Keller-Segel model has
since received much attention and there have been many contributions on this
subject such as aggregations, dynamics of blow-ups, travelling waves. See
\cite{A},\cite{B},\cite{BCKSV},\cite{HV},\cite{HP},\cite{HS},\cite{HHL},
\cite{JL},\cite{LR},\cite{MTKU},\cite{NSY},\cite{OS},\cite{P} for related
results. Linear stability and instability of stationary solutions with more
general nonlinearity was studied in \cite{S} using bifurcation analysis.
However, nonlinear evolution of the pattern formation has yet been fully
understood for the Keller-Segel model, to the authors' knowledge.

We rigorously prove that linear fastest growing modes determine unstable
patterns for the full Keller-Segel system (\ref{nonlinear1}) and
(\ref{nonlinear2}), over a time period of the order $\ln\frac{1}{\delta}.$
Each initial perturbation certainly can behaves drastically differently from
another, which gives rise to the richness of patterns. On the other hand, the
dominating linear dynamics over a fixed finite dimensional space of maximal
growing modes ensures that there is a common characteristic pattern for a
general class initial data. Therefore, we believe that our result indeed
provide a mathematical description for the pattern formation in the
Keller-Segel model.

Our paper stems from a program to study various nonlinear instabilities for
non-dissipative systems arising in mathematical physics \cite{GS},\cite{G},
\cite{BGS},\cite{HG}, where severe higher order perturbations (unbounded in
the $L^{2}$ norms, for instance) occur. Indeed, for many such systems without
dissipation, the passage from linear instability to nonlinear instability is
very delicate. If there is a dominant eigenvalue, then a bootstrap argument
was developed by Strauss and the first author to prove nonlinear instability,
for the perturbation initially along the dominant eigenfunction. The key is to
try to control the nonlinear growth of higher-order energy norm for the
perturbation by the linear growth rate, up to the time $T^{\delta}.$ Very
recently in \cite{GHS}, based upon a precise linear analysis, dynamics of
general perturbation can be characterized by the linear dynamics of fastest
growing modes for unstable Kirchhoff ellipses. This marks a beginning of a
quantitative description of instability.

Our research is inspired by the work \cite{GHS}. In the presence of
dissipation, continuum spectra are absent in bounded domain, which leads to
finite number of dominant growing modes. Moreover, natural higher-order energy
estimate now can be easily combined with the bootstrap idea to control the
nonlinear term $-\chi\nabla(u\nabla v)$ in the $L^{2}$ space. Since our method
is general, we believe that such kind of pattern formation should exist for a
wide class of systems with dissipation.

\section{Bootstrap Lemma}

We state existence of local-in-time solutions for (\ref{nonlinear1}%
)-(\ref{nonlinear2}).

\begin{lemma}
(Local existence) For $s\geq1$ $\left(  d=1\right)  $ and $s\geq2$ $\left(
d=2,3\right)  $, there exist a $T>0$ and a constant $C$ depending on
$k,\bar{U},\bar{V},D,\mu,f,\chi$ such that%
\[
\left\Vert \mathbf{w}(t)\right\Vert _{H^{s}}\leq C\left\Vert \mathbf{w}\left(
0\right)  \right\Vert _{H^{s}}.
\]

\end{lemma}

We now derive the following energy estimates for $d$-dimensional chemotaxis
model with $d=1,2,3$.

\begin{lemma}
Suppose that $[u\left(  x\mathbf{,}t\right)  ,v\left(  x,t\right)  ]$ is a
solution to the full system (\ref{nonlinear1})-(\ref{nonlinear2}). Then
\begin{align*}
&  \frac{1}{2}\frac{d}{dt}\sum_{\left\vert \partial\right\vert =2}\int_{%
{\mathbb T}^d%
}\left\{  |\partial u|^{2}+\frac{\left(  \bar{U}\chi\right)  ^{2}}{D\mu
}|\partial v|^{2}\right\}  d\mathbf{x}\\
&  +\sum_{\left\vert \partial\right\vert =2}\int_{%
{\mathbb T}^d%
}\left\{  \frac{\mu}{4}\left\vert \nabla\partial u\right\vert ^{2}%
+\frac{\left(  \bar{U}\chi\right)  ^{2}}{2\mu}\left\vert \nabla\partial
v\right\vert ^{2}\right\}  d\mathbf{x}+\frac{Ak}{2}\sum_{\left\vert
\alpha\right\vert =2}\int_{%
{\mathbb T}^d%
}|\partial v|^{2}\\
&  \leq C_{0}||\mathbf{w}||_{H^{2}}||\nabla^{3}\mathbf{w}||^{2}+C_{2}%
||u\mathbf{||}^{2}.
\end{align*}
where $C_{0}$ is the universal constant while $C_{2}=$ $\frac{\bar{U}^{6}%
\chi^{6}f^{6}}{2D^{3}\mu^{5}k^{3}}$.
\end{lemma}

\begin{proof}
We first notice that the Keller-Segel equation preserves the evenness of the
solution $\mathbf{w}(x,t),$ i.e., if $\mathbf{w}(x,t)$ is a solution, then
$\mathbf{w}(-x_{i},t)$ is also a solution. We can regard the Neumann problem
as a special case with evenness of the periodic problem by standard way of
even extension $\mathbf{w}(x,t)$ with respect to one of the $x_{i}$. For this
reason we may assume periodicity at the boundary of the extended
$2\mathbf{T}^{3}\equiv(-\pi,\pi)^{d}$. Since now there is no contributions
from the boundaries, we can take second order $\partial$-derivative of
(\ref{nonlinear1}) and add $A\times\partial$ of (\ref{nonlinear2}) to get%
\begin{align*}
&  \frac{1}{2}\frac{d}{dt}\int_{2%
{\mathbb T}^d%
}\left\{  |\partial u|^{2}+A|\partial v|^{2}\right\}  \\
&  +\int_{2%
{\mathbb T}^d%
}\{\mu\left\vert \nabla\partial u\right\vert ^{2}+AD\left\vert \nabla\partial
v\right\vert ^{2}\mathbf{-}\chi\bar{U}\nabla\partial v\nabla\partial
u\}+Ak\int_{2%
{\mathbb T}^d%
}|\partial v|^{2}\\
&  =\chi\int_{2%
{\mathbb T}^d%
}\partial\{u\nabla v\}\nabla\partial u\ +Af\int_{2%
{\mathbb T}^d%
}\partial u\partial v\\
&  \equiv I_{1}+I_{2},
\end{align*}
where the constant $A$ is given in (\ref{a}). As in (\ref{lowerbound}), the
second integrand is bounded below by
\[
\frac{\mu}{2}\left\vert \nabla\partial u\right\vert ^{2}+\frac{\left(  \bar
{U}\chi\right)  ^{2}}{2\mu}\left\vert \nabla\partial v\right\vert
^{2}\mathbf{.}%
\]
\ The nonlinear term $I_{1}$ is bounded by
\begin{align*}
I_{1} &  =\int\left\vert \partial\left(  u\nabla v\right)  \cdot\nabla\partial
u\right\vert d\mathbf{x}\\
\leq &  \left\Vert u\right\Vert _{L^{\infty}}\left\Vert \nabla\partial
v\right\Vert \left\Vert \nabla\partial u\right\Vert +\left\Vert \nabla
u\right\Vert _{L^{\infty}}\left\Vert \partial v\right\Vert \left\Vert
\nabla\partial u\right\Vert \\
&  +\left\Vert \partial u\right\Vert \left\Vert \nabla v\right\Vert
_{L^{\infty}}\left\Vert \nabla\partial u\right\Vert .
\end{align*}
We apply the following the Sobolev imbedding to control $\left\Vert
u\right\Vert _{L^{\infty}}$
\begin{equation}
\left\Vert g\right\Vert _{L^{\infty}\left(  2%
{\mathbb T}^d%
\right)  }\leq C_{0}\left\Vert g\right\Vert _{H^{2}\left(  2%
{\mathbb T}^d%
\right)  },\label{Sobolev1.5}%
\end{equation}
for $d\leq3.$ Moreover, from the periodic boundary conditions,
\[
\int_{2%
{\mathbb T}^d%
}\nabla u=\int_{2%
{\mathbb T}^d%
}\nabla v=0,
\]
we also use the Poincare inequality%
\begin{equation}
||g||\leq\left\Vert g\right\Vert _{L^{4}(2%
{\mathbb T}^d%
)}\leq C_{0}\left\Vert \nabla g\right\Vert \text{ \ \ \ if }d\leq
3,\label{Sobolev2}%
\end{equation}
to further get%
\[
||\nabla u||_{L^{\infty}}+||\nabla v||_{L^{\infty}}\leq C_{0}\{\left\Vert
\nabla u\right\Vert _{H^{2}}+\left\Vert \nabla v\right\Vert _{H^{2}}\}\leq
C_{0}\sum_{|\partial|=2}||\partial\nabla\mathbf{w}||,
\]
where $C_{0}$ is a universal constant. Hence $I_{1}\leq C_{0}||\mathbf{w}%
||_{H^{2}}||\nabla^{3}\mathbf{w}||^{2}$ as desired.

$\ \ \ \ \ \ $Finally, $I_{2}$ is simply bounded by%
\[
I_{2}=Af\int\partial u\partial v\leq\frac{Af^{2}}{2k}||\partial u||^{2}%
+\frac{Ak}{2}||\partial v||^{2}%
\]
By the interpolation between $\left\Vert \nabla\partial u\right\Vert $ and
$||u\mathbf{||}$ , the first term above is bounded by
\[
\frac{Af^{2}}{2k}\{a\left\Vert \nabla\partial u\right\Vert ^{2}+\frac
{1}{4a^{2}}||u\mathbf{||}^{2}\}
\]
for any $a>0.$ We can choose $a$ such that $\frac{Af^{2}}{2k}a=\frac{1}{4}\mu
$. Collecting terms, we conclude the proof.
\end{proof}

We are now ready to establish the bootstrap lemma, which controls the $H^{2}$
growth of $\mathbf{w}(x,t)$ in term of its $L^{2}$ growth$.$

\begin{lemma}
Suppose that $\mathbf{w}(x,t)$ is a solution to the full system
(\ref{nonlinear1})-(\ref{nonlinear2}) such that for $0\leq t\leq T,$
\[
||\mathbf{w}(\cdot,t)||_{H^{2}}\leq\frac{1}{C_{0}}\min\left\{  \frac{\mu}%
{4},\frac{\left(  \bar{U}\chi\right)  ^{2}}{2\mu}\right\}
\]
and
\begin{equation}
||\mathbf{w}(\cdot,t)||\leq2C_{1}e^{\lambda_{\max}t}||\mathbf{w}(\cdot,0)||,
\label{l2growth}%
\end{equation}
then we have for $0\leq t\leq T,$
\[
||\mathbf{w}(t)||_{H^{2}}^{2}\leq C_{3}\{||\mathbf{w}(0)||_{H^{2}}%
^{2}+e^{2\lambda_{\max}t}||\mathbf{w}(\cdot,0)||^{2}\}
\]
where $C_{3}$ $=C_{1}^{2}\max\{\frac{\left(  \bar{U}\chi\right)  ^{2}}{D\mu
},\frac{D\mu}{\left(  \bar{U}\chi\right)  ^{2}}\}\times\max\{\frac{4C_{2}%
}{\lambda_{\text{max}}},1\}\geq1.$
\end{lemma}

\begin{proof}
It suffices to only consider the second-order derivatives of $\mathbf{w}%
(x,t).$ From the previous lemma and our assumption for $||\mathbf{w}||_{H^{2}%
},$ we deduce that for $0\leq t\leq T$
\[
\frac{1}{2}\frac{d}{dt}\sum_{\left\vert \partial\right\vert =2}\int_{%
{\mathbb T}^d%
}\left\{  |\partial u|^{2}+\frac{\left(  \bar{U}\chi\right)  ^{2}}{D\mu
}|\partial v|^{2}\right\}  d\mathbf{x}\leq C_{2}||u\mathbf{||}^{2}.
\]
So that by (\ref{l2growth}) and an integration from $0$ to $t,$ we have
\begin{align*}
&  \sum_{|\partial|=2}\int_{%
{\mathbb T}^d%
}\left\{  |\partial u(t)|^{2}+\frac{\left(  \bar{U}\chi\right)  ^{2}}{D\mu
}|\partial v(t)|^{2}\right\}  \\
&  \leq\sum_{|\partial|=2}\int_{%
{\mathbb T}^d%
}\left\{  |\partial u(0)|^{2}+\frac{\left(  \bar{U}\chi\right)  ^{2}}{D\mu
}|\partial v(0)|^{2}\right\}  +\frac{4C_{2}C_{1}^{2}}{\lambda_{\text{max}}%
}e^{2\lambda_{\max}t}||\mathbf{w}(\cdot,0)||^{2},
\end{align*}
for $0\leq t\leq T.$ Now our lemma follows directly by separating the cases of
$\frac{\left(  \bar{U}\chi\right)  ^{2}}{D\mu}\geq1$ and $\frac{\left(
\bar{U}\chi\right)  ^{2}}{D\mu}<1.$
\end{proof}

\section{Nonlinear instability and pattern formation}

We now prove our main Theorem 1:

\begin{proof}
Let $\mathbf{w}^{\delta}\left(  x,t\right)  $ be the family of solutions to
the Keller-Segel system (\ref{nonlinear1})-(\ref{nonlinear2}) with initial
data $\mathbf{w}^{\delta}\left(  x,0\right)  =\delta\mathbf{w}_{0}$. Define
$T^{\ast}$ by%
\[
T^{\ast}=\sup\left\{  t\ |~\left\Vert \mathbf{w}^{\delta}(t)-\delta e^{%
{\mathcal L}%
t}\mathbf{w}_{0}\right\Vert \leq\frac{C_{1}}{2}\delta\exp\left(  \lambda
_{\max}t\right)  \right\}  .
\]
Note that $T^{\ast}$ is well defined. We also define%
\[
T^{\ast\ast}=\sup\left\{  t\ |~\left\Vert \mathbf{w}^{\delta}(t)\right\Vert
_{H^{2}}\leq\frac{1}{C_{0}}\min\left\{  \frac{\mu}{4},\frac{\left(  \bar
{U}\chi\right)  ^{2}}{2\mu}\right\}  \right\}  .
\]
We recall $T^{\delta}$ in (\ref{theta}) where $\theta$ is chosen such that
\begin{equation}
C_{0}C_{3}\theta<\min\left\{  \frac{\lambda_{\text{max}}}{4},\frac{\mu}%
{8},\frac{\left(  \bar{U}\chi\right)  ^{2}}{4\mu}\right\}  , \label{th1}%
\end{equation}
\ 

\ \ \ \ \ We now derive estimates for $H^{2}$ norm of $\mathbf{w}^{\delta
}(x,t)$ for $0\leq t\leq\min\{T^{\ast},T^{\delta},T^{\ast\ast}\}.$ First of
all, by the definition of $T^{\ast},$ for $t\leq T^{\ast}$ and Lemma
\ref{lineargrowth}
\[
\left\Vert \mathbf{w}^{\delta}(t)\right\Vert \leq\frac{3C_{1}}{2}\delta
\exp\left(  \lambda_{\max}t\right)  .
\]
\ \ \ \ \ \ \ 

Moreover, using Lemma 4 and applying a bootstrap argument yields%
\begin{equation}
\left\Vert \mathbf{w}^{\delta}(t)\right\Vert _{H^{2}}\leq\sqrt{C_{3}}%
\{\delta||\mathbf{w}_{0}||_{H^{2}}+\delta e^{\lambda_{\max}t}\}. \label{h2}%
\end{equation}
\ 

We now establish a sharper $L^{2}$ estimate for $\mathbf{w}^{\delta}(x,t)$,
for $0\leq t\leq\min\{T^{\ast\ast},T^{\delta},T^{\ast}\}.$ We first apply
Duhamel's principle to obtain%
\[
\mathbf{w}^{\delta}\left(  t\right)  =\delta e^{%
{\mathcal L}%
t}\mathbf{w}_{0}-\int_{0}^{t}e^{%
{\mathcal L}%
\left(  t-\tau\right)  }[\nabla\cdot\left(  u^{\delta}\left(  \tau\right)
\nabla v^{\delta}\left(  \tau\right)  \right)  ,0]d\tau,
\]
Using Lemma 1, (\ref{Sobolev1.5}), (\ref{Sobolev2}), and Lemma 4 yields, for
$0\leq t\leq\min\{T^{\delta},T^{\ast\ast},T^{\ast}\}$%
\begin{align*}
&  \left\Vert \mathbf{w}^{\delta}\left(  t\right)  -\delta e^{%
{\mathcal L}%
t}\mathbf{w}_{0}\right\Vert \\
\leq &  C_{1}\int_{0}^{t}e^{\lambda_{\max}\left(  t-\tau\right)  }\left\Vert
\nabla\cdot\left(  u^{\delta}\left(  \tau\right)  \nabla v^{\delta}\left(
\tau\right)  \right)  \right\Vert d\tau\\
\leq &  C_{1}\int_{0}^{t}e^{\lambda_{\max}\left(  t-\tau\right)  }\left\Vert
u^{\delta}\left(  \tau\right)  \right\Vert _{L^{\infty}}\left\Vert \nabla
^{2}v^{\delta}\left(  \tau\right)  \right\Vert d\tau\\
&  +C_{1}\int_{0}^{t}e^{\lambda_{\max}\left(  t-\tau\right)  }\left\Vert
\nabla u^{\delta}\left(  \tau\right)  \right\Vert _{L^{4}}\left\Vert \nabla
v^{\delta}\left(  \tau\right)  \right\Vert _{L^{4}}d\tau\\
\leq &  C_{1}C_{0}\int_{0}^{t}e^{\lambda_{\max}\left(  t-\tau\right)
}\left\Vert \mathbf{w}^{\delta}\left(  \tau\right)  \right\Vert _{H^{2}}%
^{2}d\tau.
\end{align*}
By our choice of $t\leq\min\{T^{\ast},T^{\ast\ast},T^{\delta}\},$ it is
further bounded by%
\begin{align}
&  \left\Vert \mathbf{w}^{\delta}\left(  t\right)  -\delta e^{%
{\mathcal L}%
t}\mathbf{w}_{0}\right\Vert \label{l2}\\
&  \leq C_{1}C_{0}C_{3}\int_{0}^{t}e^{\lambda_{\max}\left(  t-\tau\right)
}\{\delta^{2}||\mathbf{w}_{0}||_{H^{2}}^{2}+\delta^{2}e^{2\lambda_{\max}\tau
}\}d\tau\nonumber\\
&  \leq C_{1}C_{0}C_{3}\{\frac{||\mathbf{w}_{0}||_{H^{2}}^{2}\delta}%
{\lambda_{\text{max}}}+\frac{1}{\lambda_{\text{max}}}\delta e^{\lambda_{\max
}t}\}\delta e^{\lambda_{\max}t}.\nonumber
\end{align}
\ \ \ \ \ \ 

\ \ We now prove by contradiction that for $\delta$ sufficiently small,
\[
T^{\delta}=\min\{T^{\delta},T^{\ast},T^{\ast\ast}\},
\]
and therefore our theorem follows by further separating $\mathbf{q}\in
\Omega_{\max}$ and move $\mathbf{q}\notin\Omega_{\max}$ in (\ref{l}) to the
right hand side .

\ \ \ If $T^{\ast\ast}$ is the smallest, we can let $t=T^{\ast\ast}\leq
T^{\delta}$ in (\ref{h2})
\begin{align*}
\left\Vert \mathbf{w}^{\delta}(T^{\ast\ast})\right\Vert _{H^{2}}  &  \leq
\sqrt{C_{3}}\{\delta||\mathbf{w}_{0}||_{H^{2}}+\delta e^{\lambda_{\max
}T^{\delta}}\}\\
&  =\sqrt{C_{3}}\{\delta||\mathbf{w}_{0}||_{H^{2}}+\theta\}\\
&  <\frac{1}{C_{0}}\min\{\frac{\mu}{4},\frac{\left(  \bar{U}\chi\right)  ^{2}%
}{2\mu}\},
\end{align*}
for $\sqrt{C_{3}}\delta||\mathbf{w}_{0}||_{H^{2}}\leq\frac{1}{2C_{0}}%
\min\{\frac{\mu}{4},\frac{\left(  \bar{U}\chi\right)  ^{2}}{2\mu}\},$ by our
choice of $\theta$ in (\ref{th1}) with $C_{3}\geq1$. This is a contradiction
to the definition of $T^{\ast\ast}.$

\ \ \ On the other hand, if $T^{\ast}$ is the smallest, we let $t=T^{\ast}$ in
(\ref{l2}) to get%
\begin{align*}
&  \left\Vert \mathbf{w}^{\delta}\left(  T^{\ast}\right)  -\delta e^{%
{\mathcal L}%
T^{\ast}}\mathbf{w}_{0}\right\Vert \\
&  \leq C_{1}C_{0}C_{3}\{\frac{||\mathbf{w}_{0}||_{H^{2}}^{2}\delta}%
{\lambda_{\text{max}}}+\frac{1}{\lambda_{\text{max}}}\delta e^{\lambda_{\max
}T^{\delta}}\}\delta e^{\lambda_{\max}T^{\ast}}\\
&  \leq C_{1}C_{0}C_{3}\{\frac{||\mathbf{w}_{0}||_{H^{2}}^{2}\delta}%
{\lambda_{\text{max}}}+\frac{\theta}{\lambda_{\text{max}}}\}\delta
e^{\lambda_{\max}T^{\ast}}\\
&  <\frac{C_{1}}{2}\delta e^{\lambda_{\max}T^{\ast}},
\end{align*}
for $C_{0}C_{3}\frac{||\mathbf{w}_{0}||_{H^{2}}^{2}\delta}{\lambda
_{\text{max}}}<1/4$ for $\delta$ small, by our choice of $\theta$ in
(\ref{th1}). This again contradicts the definition of $T^{\ast}$ and our
theorem follows.
\end{proof}

\end{document}